\documentclass[11pt]{article}

    \usepackage[breakable]{tcolorbox}
    \usepackage{parskip} 
    
    \usepackage{iftex}
    \ifPDFTeX
    	\usepackage[T1]{fontenc}
    	\usepackage{mathpazo}
    \else
    	\usepackage{fontspec}
    \fi

    \usepackage{graphicx}
    
    \usepackage{caption}
    \DeclareCaptionFormat{nocaption}{}
    \usepackage[Export]{adjustbox} 
    \adjustboxset{max size={0.9\linewidth}{0.9\paperheight}}
    \usepackage{float}
    \usepackage{xcolor} 
    \usepackage{enumerate} 
    \usepackage{geometry} 
    \usepackage{amsmath} 
    \usepackage{amssymb} 
    \usepackage{textcomp} 
    \AtBeginDocument{%
        \def\PYZsq{\textquotesingle}
    }
    \usepackage{upquote} 
    \usepackage{eurosym} 
    \usepackage[mathletters]{ucs} 
    \usepackage{fancyvrb} 
    \usepackage{grffile} 
    \makeatletter 
    \def\Gread@@xetex#1{%
      \IfFileExists{"\Gin@base".bb}%
      {\Gread@eps{\Gin@base.bb}}%
      {\Gread@@xetex@aux#1}%
    }
    \makeatother

    \usepackage{hyperref}
    \usepackage{titling}
    \usepackage{longtable} 
    \usepackage{booktabs}  
    \usepackage[inline]{enumitem} 
    \usepackage[normalem]{ulem} 
    \usepackage{mathrsfs}

    \definecolor{urlcolor}{rgb}{0,.145,.698}
    \definecolor{linkcolor}{rgb}{.71,0.21,0.01}
    \definecolor{citecolor}{rgb}{.12,.54,.11}

    \definecolor{ansi-black}{HTML}{3E424D}
    \definecolor{ansi-black-intense}{HTML}{282C36}
    \definecolor{ansi-red}{HTML}{E75C58}
    \definecolor{ansi-red-intense}{HTML}{B22B31}
    \definecolor{ansi-green}{HTML}{00A250}
    \definecolor{ansi-green-intense}{HTML}{007427}
    \definecolor{ansi-yellow}{HTML}{DDB62B}
    \definecolor{ansi-yellow-intense}{HTML}{B27D12}
    \definecolor{ansi-blue}{HTML}{208FFB}
    \definecolor{ansi-blue-intense}{HTML}{0065CA}
    \definecolor{ansi-magenta}{HTML}{D160C4}
    \definecolor{ansi-magenta-intense}{HTML}{A03196}
    \definecolor{ansi-cyan}{HTML}{60C6C8}
    \definecolor{ansi-cyan-intense}{HTML}{258F8F}
    \definecolor{ansi-white}{HTML}{C5C1B4}
    \definecolor{ansi-white-intense}{HTML}{A1A6B2}
    \definecolor{ansi-default-inverse-fg}{HTML}{FFFFFF}
    \definecolor{ansi-default-inverse-bg}{HTML}{000000}

    \DefineVerbatimEnvironment{Highlighting}{Verbatim}{commandchars=\\\{\}}

    \let\Oldtex\TeX
    \let\Oldlatex\LaTeX
    \renewcommand{\TeX}{\textrm{\Oldtex}}
    \renewcommand{\LaTeX}{\textrm{\Oldlatex}}
    
\title{Computing Eigenvectors from Eigenvalues In an Arbitrary Orthonormal Basis}

\date{19 November 2019}

\author{John Lakness}

\makeatletter
\def\PY@reset{\let\PY@it=\relax \let\PY@bf=\relax%
    \let\PY@ul=\relax \let\PY@tc=\relax%
    \let\PY@bc=\relax \let\PY@ff=\relax}
\def\PY@tok#1{\csname PY@tok@#1\endcsname}
\def\PY@toks#1+{\ifx\relax#1\empty\else%
    \PY@tok{#1}\expandafter\PY@toks\fi}
\def\PY@do#1{\PY@bc{\PY@tc{\PY@ul{%
    \PY@it{\PY@bf{\PY@ff{#1}}}}}}}
\def\PY#1#2{\PY@reset\PY@toks#1+\relax+\PY@do{#2}}

\expandafter\def\csname PY@tok@w\endcsname{\def\PY@tc##1{\textcolor[rgb]{0.73,0.73,0.73}{##1}}}
\expandafter\def\csname PY@tok@c\endcsname{\let\PY@it=\textit\def\PY@tc##1{\textcolor[rgb]{0.25,0.50,0.50}{##1}}}
\expandafter\def\csname PY@tok@cp\endcsname{\def\PY@tc##1{\textcolor[rgb]{0.74,0.48,0.00}{##1}}}
\expandafter\def\csname PY@tok@k\endcsname{\let\PY@bf=\textbf\def\PY@tc##1{\textcolor[rgb]{0.00,0.50,0.00}{##1}}}
\expandafter\def\csname PY@tok@kp\endcsname{\def\PY@tc##1{\textcolor[rgb]{0.00,0.50,0.00}{##1}}}
\expandafter\def\csname PY@tok@kt\endcsname{\def\PY@tc##1{\textcolor[rgb]{0.69,0.00,0.25}{##1}}}
\expandafter\def\csname PY@tok@o\endcsname{\def\PY@tc##1{\textcolor[rgb]{0.40,0.40,0.40}{##1}}}
\expandafter\def\csname PY@tok@ow\endcsname{\let\PY@bf=\textbf\def\PY@tc##1{\textcolor[rgb]{0.67,0.13,1.00}{##1}}}
\expandafter\def\csname PY@tok@nb\endcsname{\def\PY@tc##1{\textcolor[rgb]{0.00,0.50,0.00}{##1}}}
\expandafter\def\csname PY@tok@nf\endcsname{\def\PY@tc##1{\textcolor[rgb]{0.00,0.00,1.00}{##1}}}
\expandafter\def\csname PY@tok@nc\endcsname{\let\PY@bf=\textbf\def\PY@tc##1{\textcolor[rgb]{0.00,0.00,1.00}{##1}}}
\expandafter\def\csname PY@tok@nn\endcsname{\let\PY@bf=\textbf\def\PY@tc##1{\textcolor[rgb]{0.00,0.00,1.00}{##1}}}
\expandafter\def\csname PY@tok@ne\endcsname{\let\PY@bf=\textbf\def\PY@tc##1{\textcolor[rgb]{0.82,0.25,0.23}{##1}}}
\expandafter\def\csname PY@tok@nv\endcsname{\def\PY@tc##1{\textcolor[rgb]{0.10,0.09,0.49}{##1}}}
\expandafter\def\csname PY@tok@no\endcsname{\def\PY@tc##1{\textcolor[rgb]{0.53,0.00,0.00}{##1}}}
\expandafter\def\csname PY@tok@nl\endcsname{\def\PY@tc##1{\textcolor[rgb]{0.63,0.63,0.00}{##1}}}
\expandafter\def\csname PY@tok@ni\endcsname{\let\PY@bf=\textbf\def\PY@tc##1{\textcolor[rgb]{0.60,0.60,0.60}{##1}}}
\expandafter\def\csname PY@tok@na\endcsname{\def\PY@tc##1{\textcolor[rgb]{0.49,0.56,0.16}{##1}}}
\expandafter\def\csname PY@tok@nt\endcsname{\let\PY@bf=\textbf\def\PY@tc##1{\textcolor[rgb]{0.00,0.50,0.00}{##1}}}
\expandafter\def\csname PY@tok@nd\endcsname{\def\PY@tc##1{\textcolor[rgb]{0.67,0.13,1.00}{##1}}}
\expandafter\def\csname PY@tok@s\endcsname{\def\PY@tc##1{\textcolor[rgb]{0.73,0.13,0.13}{##1}}}
\expandafter\def\csname PY@tok@sd\endcsname{\let\PY@it=\textit\def\PY@tc##1{\textcolor[rgb]{0.73,0.13,0.13}{##1}}}
\expandafter\def\csname PY@tok@si\endcsname{\let\PY@bf=\textbf\def\PY@tc##1{\textcolor[rgb]{0.73,0.40,0.53}{##1}}}
\expandafter\def\csname PY@tok@se\endcsname{\let\PY@bf=\textbf\def\PY@tc##1{\textcolor[rgb]{0.73,0.40,0.13}{##1}}}
\expandafter\def\csname PY@tok@sr\endcsname{\def\PY@tc##1{\textcolor[rgb]{0.73,0.40,0.53}{##1}}}
\expandafter\def\csname PY@tok@ss\endcsname{\def\PY@tc##1{\textcolor[rgb]{0.10,0.09,0.49}{##1}}}
\expandafter\def\csname PY@tok@sx\endcsname{\def\PY@tc##1{\textcolor[rgb]{0.00,0.50,0.00}{##1}}}
\expandafter\def\csname PY@tok@m\endcsname{\def\PY@tc##1{\textcolor[rgb]{0.40,0.40,0.40}{##1}}}
\expandafter\def\csname PY@tok@gh\endcsname{\let\PY@bf=\textbf\def\PY@tc##1{\textcolor[rgb]{0.00,0.00,0.50}{##1}}}
\expandafter\def\csname PY@tok@gu\endcsname{\let\PY@bf=\textbf\def\PY@tc##1{\textcolor[rgb]{0.50,0.00,0.50}{##1}}}
\expandafter\def\csname PY@tok@gd\endcsname{\def\PY@tc##1{\textcolor[rgb]{0.63,0.00,0.00}{##1}}}
\expandafter\def\csname PY@tok@gi\endcsname{\def\PY@tc##1{\textcolor[rgb]{0.00,0.63,0.00}{##1}}}
\expandafter\def\csname PY@tok@gr\endcsname{\def\PY@tc##1{\textcolor[rgb]{1.00,0.00,0.00}{##1}}}
\expandafter\def\csname PY@tok@ge\endcsname{\let\PY@it=\textit}
\expandafter\def\csname PY@tok@gs\endcsname{\let\PY@bf=\textbf}
\expandafter\def\csname PY@tok@gp\endcsname{\let\PY@bf=\textbf\def\PY@tc##1{\textcolor[rgb]{0.00,0.00,0.50}{##1}}}
\expandafter\def\csname PY@tok@go\endcsname{\def\PY@tc##1{\textcolor[rgb]{0.53,0.53,0.53}{##1}}}
\expandafter\def\csname PY@tok@gt\endcsname{\def\PY@tc##1{\textcolor[rgb]{0.00,0.27,0.87}{##1}}}
\expandafter\def\csname PY@tok@err\endcsname{\def\PY@bc##1{\setlength{\fboxsep}{0pt}\fcolorbox[rgb]{1.00,0.00,0.00}{1,1,1}{\strut ##1}}}
\expandafter\def\csname PY@tok@kc\endcsname{\let\PY@bf=\textbf\def\PY@tc##1{\textcolor[rgb]{0.00,0.50,0.00}{##1}}}
\expandafter\def\csname PY@tok@kd\endcsname{\let\PY@bf=\textbf\def\PY@tc##1{\textcolor[rgb]{0.00,0.50,0.00}{##1}}}
\expandafter\def\csname PY@tok@kn\endcsname{\let\PY@bf=\textbf\def\PY@tc##1{\textcolor[rgb]{0.00,0.50,0.00}{##1}}}
\expandafter\def\csname PY@tok@kr\endcsname{\let\PY@bf=\textbf\def\PY@tc##1{\textcolor[rgb]{0.00,0.50,0.00}{##1}}}
\expandafter\def\csname PY@tok@bp\endcsname{\def\PY@tc##1{\textcolor[rgb]{0.00,0.50,0.00}{##1}}}
\expandafter\def\csname PY@tok@fm\endcsname{\def\PY@tc##1{\textcolor[rgb]{0.00,0.00,1.00}{##1}}}
\expandafter\def\csname PY@tok@vc\endcsname{\def\PY@tc##1{\textcolor[rgb]{0.10,0.09,0.49}{##1}}}
\expandafter\def\csname PY@tok@vg\endcsname{\def\PY@tc##1{\textcolor[rgb]{0.10,0.09,0.49}{##1}}}
\expandafter\def\csname PY@tok@vi\endcsname{\def\PY@tc##1{\textcolor[rgb]{0.10,0.09,0.49}{##1}}}
\expandafter\def\csname PY@tok@vm\endcsname{\def\PY@tc##1{\textcolor[rgb]{0.10,0.09,0.49}{##1}}}
\expandafter\def\csname PY@tok@sa\endcsname{\def\PY@tc##1{\textcolor[rgb]{0.73,0.13,0.13}{##1}}}
\expandafter\def\csname PY@tok@sb\endcsname{\def\PY@tc##1{\textcolor[rgb]{0.73,0.13,0.13}{##1}}}
\expandafter\def\csname PY@tok@sc\endcsname{\def\PY@tc##1{\textcolor[rgb]{0.73,0.13,0.13}{##1}}}
\expandafter\def\csname PY@tok@dl\endcsname{\def\PY@tc##1{\textcolor[rgb]{0.73,0.13,0.13}{##1}}}
\expandafter\def\csname PY@tok@s2\endcsname{\def\PY@tc##1{\textcolor[rgb]{0.73,0.13,0.13}{##1}}}
\expandafter\def\csname PY@tok@sh\endcsname{\def\PY@tc##1{\textcolor[rgb]{0.73,0.13,0.13}{##1}}}
\expandafter\def\csname PY@tok@s1\endcsname{\def\PY@tc##1{\textcolor[rgb]{0.73,0.13,0.13}{##1}}}
\expandafter\def\csname PY@tok@mb\endcsname{\def\PY@tc##1{\textcolor[rgb]{0.40,0.40,0.40}{##1}}}
\expandafter\def\csname PY@tok@mf\endcsname{\def\PY@tc##1{\textcolor[rgb]{0.40,0.40,0.40}{##1}}}
\expandafter\def\csname PY@tok@mh\endcsname{\def\PY@tc##1{\textcolor[rgb]{0.40,0.40,0.40}{##1}}}
\expandafter\def\csname PY@tok@mi\endcsname{\def\PY@tc##1{\textcolor[rgb]{0.40,0.40,0.40}{##1}}}
\expandafter\def\csname PY@tok@il\endcsname{\def\PY@tc##1{\textcolor[rgb]{0.40,0.40,0.40}{##1}}}
\expandafter\def\csname PY@tok@mo\endcsname{\def\PY@tc##1{\textcolor[rgb]{0.40,0.40,0.40}{##1}}}
\expandafter\def\csname PY@tok@ch\endcsname{\let\PY@it=\textit\def\PY@tc##1{\textcolor[rgb]{0.25,0.50,0.50}{##1}}}
\expandafter\def\csname PY@tok@cm\endcsname{\let\PY@it=\textit\def\PY@tc##1{\textcolor[rgb]{0.25,0.50,0.50}{##1}}}
\expandafter\def\csname PY@tok@cpf\endcsname{\let\PY@it=\textit\def\PY@tc##1{\textcolor[rgb]{0.25,0.50,0.50}{##1}}}
\expandafter\def\csname PY@tok@c1\endcsname{\let\PY@it=\textit\def\PY@tc##1{\textcolor[rgb]{0.25,0.50,0.50}{##1}}}
\expandafter\def\csname PY@tok@cs\endcsname{\let\PY@it=\textit\def\PY@tc##1{\textcolor[rgb]{0.25,0.50,0.50}{##1}}}

\def\PYZus{\char`\_}

\def\PYZlt{\char`\<}

\def\PYZhy{\char`\-}
\def\PYZsq{\char`\'}

\makeatother

    \makeatletter
        \newbox\Wrappedcontinuationbox 
        \newbox\Wrappedvisiblespacebox 
        \newcommand*\Wrappedvisiblespace {\textcolor{red}{\textvisiblespace}} 
        \newcommand*\Wrappedcontinuationsymbol {\textcolor{red}{\llap{\tiny$\m@th\hookrightarrow$}}} 
        \newcommand*\Wrappedcontinuationindent {3ex } 
        \newcommand*\Wrappedafterbreak {\kern\Wrappedcontinuationindent\copy\Wrappedcontinuationbox} 
        \newcommand*\Wrappedbreaksatspecials {%
            \def\PYGZus{\discretionary{\char`\_}{\Wrappedafterbreak}{\char`\_}}%
            \def\PYGZob{\discretionary{}{\Wrappedafterbreak\char`\{}{\char`\{}}%
            \def\PYGZcb{\discretionary{\char`\}}{\Wrappedafterbreak}{\char`\}}}%
            \def\PYGZca{\discretionary{\char`\^}{\Wrappedafterbreak}{\char`\^}}%
            \def\PYGZam{\discretionary{\char`\&}{\Wrappedafterbreak}{\char`\&}}%
            \def\PYGZlt{\discretionary{}{\Wrappedafterbreak\char`\<}{\char`\<}}%
            \def\PYGZgt{\discretionary{\char`\>}{\Wrappedafterbreak}{\char`\>}}%
            \def\PYGZsh{\discretionary{}{\Wrappedafterbreak\char`\#}{\char`\#}}%
            \def\PYGZpc{\discretionary{}{\Wrappedafterbreak\char`\%}{\char`\%}}%
            \def\PYGZdl{\discretionary{}{\Wrappedafterbreak\char`\$}{\char`\$}}%
            \def\PYGZhy{\discretionary{\char`\-}{\Wrappedafterbreak}{\char`\-}}%
            \def\PYGZsq{\discretionary{}{\Wrappedafterbreak\textquotesingle}{\textquotesingle}}%
            \def\PYGZdq{\discretionary{}{\Wrappedafterbreak\char`\"}{\char`\"}}%
            \def\PYGZti{\discretionary{\char`\~}{\Wrappedafterbreak}{\char`\~}}%
        } 
        \newcommand*\Wrappedbreaksatpunct {%
            \lccode`\~`\.\lowercase{\def~}{\discretionary{\hbox{\char`\.}}{\Wrappedafterbreak}{\hbox{\char`\.}}}%
            \lccode`\~`\,\lowercase{\def~}{\discretionary{\hbox{\char`\,}}{\Wrappedafterbreak}{\hbox{\char`\,}}}%
            \lccode`\~`\;\lowercase{\def~}{\discretionary{\hbox{\char`\;}}{\Wrappedafterbreak}{\hbox{\char`\;}}}%
            \lccode`\~`\:\lowercase{\def~}{\discretionary{\hbox{\char`\:}}{\Wrappedafterbreak}{\hbox{\char`\:}}}%
            \lccode`\~`\?\lowercase{\def~}{\discretionary{\hbox{\char`\?}}{\Wrappedafterbreak}{\hbox{\char`\?}}}%
            \lccode`\~`\!\lowercase{\def~}{\discretionary{\hbox{\char`\!}}{\Wrappedafterbreak}{\hbox{\char`\!}}}%
            \lccode`\~`\/\lowercase{\def~}{\discretionary{\hbox{\char`\/}}{\Wrappedafterbreak}{\hbox{\char`\/}}}%
            \catcode`\.\active
            \catcode`\,\active 
            \catcode`\;\active
            \catcode`\:\active
            \catcode`\?\active
            \catcode`\!\active
            \catcode`\/\active 
            \lccode`\~`\~ 	
        }
    \makeatother

    \let\OriginalVerbatim=\Verbatim
    \makeatletter
    \renewcommand{\Verbatim}[1][1]{%
        \sbox\Wrappedcontinuationbox {\Wrappedcontinuationsymbol}%
        \sbox\Wrappedvisiblespacebox {\FV@SetupFont\Wrappedvisiblespace}%
        \def\FancyVerbFormatLine ##1{\hsize\linewidth
            \vtop{\raggedright\hyphenpenalty\z@\exhyphenpenalty\z@
                \doublehyphendemerits\z@\finalhyphendemerits\z@
                \strut ##1\strut}%
        }%
        \def\FV@Space {%
            \nobreak\hskip\z@ plus\fontdimen3\font minus\fontdimen4\font
            \discretionary{\copy\Wrappedvisiblespacebox}{\Wrappedafterbreak}
            {\kern\fontdimen2\font}%
        }%
        
        \Wrappedbreaksatspecials
        \OriginalVerbatim[#1,codes*=\Wrappedbreaksatpunct]%
    }
    \makeatother

    \definecolor{incolor}{HTML}{303F9F}
    \definecolor{outcolor}{HTML}{D84315}
    \definecolor{cellborder}{HTML}{CFCFCF}
    \definecolor{cellbackground}{HTML}{F7F7F7}
    
    \makeatletter
    \newcommand{\boxspacing}{\kern\kvtcb@left@rule\kern\kvtcb@boxsep}
    \makeatother
    \newcommand{\prompt}[4]{
        \ttfamily\llap{{\color{#2}[#3]:\hspace{3pt}#4}}\vspace{-\baselineskip}
    }

    \hypersetup{
      breaklinks=true,  
      colorlinks=true,
      urlcolor=urlcolor,
      linkcolor=linkcolor,
      citecolor=citecolor,
      }
    
\begin{document}
    
    \maketitle

    \hypertarget{abstract}{%
\subsection{Abstract}\label{abstract}}

The method of computing eigenvectors from eigenvalues of submatrices can
be shown as equivalent to a method of computing the constraint which
achieves specified stationary values of a quadratic optimization.
Similarly, we show computation of eigenvectors of an orthonormal basis
projection using eigenvalues of sub-projections.

    \hypertarget{eigenvector-element-magnitude-from-sub-matrix-eigenvalues}{%
\section{Eigenvector Element Magnitude from Sub-Matrix
Eigenvalues}\label{eigenvector-element-magnitude-from-sub-matrix-eigenvalues}}

A recent result proposed by Denton and proven by Tao showed
\cite{denton2019eigenvectors}:

Let \(A\) be a Hermitian matrix of size \(n \times n\) where
\(A = QWQ^*\) is its eigendecomposition having eigenvectors \(q_i\) with
\(j^{th}\) element \(q_{ij}\) and corresponding eigenvalues \(w_i\). If
\(M_{j}\) is the \((n-1) \times (n-1)\) matrix formed by deleting the
\(j^{th}\) row and column of \(A\) where \(M_j=Q_jX_jQ_j^*\) with
eigenvalues \(x_{jk}\), then:
\[ q_{ij}^*q_{ij} = \frac{\prod_{k=1;k\neq i}^{n}{w_i-w_k}}{\prod_{k=1}^{n-1}{w_i-x_{jk}}} \]

    \hypertarget{quadratic-program-constraints-from-stationary-values}{%
\section{Quadratic Program Constraints from Stationary
Values}\label{quadratic-program-constraints-from-stationary-values}}

Previously results by Golub had shown \cite{10.2307/2028604}:

\begin{enumerate}
\def\labelenumi{\arabic{enumi}.}
\item
  Let \(A\) be a real symmetric matrix of order \(n\). Let \(c\) an
  n-vector with \(c^Tc=1\). The stationary values of \(x^TAx\) subject
  to \(x^Tx=1\) and \(c^Tx=0\) are the eigenvalues of \(PAP\) where
  \(P=I-cc^T\).
\item
  If \(A\) has eigendecomposition \(A = QWQ^T\) with eigenvalues
  \(w_i\), \(PAP\) has non-zero eigenvalues of \(x_k\), and and \(Qd=c\)
  with \(d_j\) the \(j^{th}\) element of \(d\),then
  \[ d_j^2= \frac{\prod_{k=1;k\neq j}^{n}{w_j-w_k}}{\prod_{k=1}^{n-1}{w_j-x_{k}}} \]
\end{enumerate}

This allows us to find the constraint vector \(c\) which produces a set
of arbitrary stationary values \(x_k\) by the relationship \(Qd=c\)
where there are two possibilities for each \(d_j\) from our computation
of \(d_j^2\).

    \hypertarget{equivalence}{%
\section{Equivalence}\label{equivalence}}

We can show that these two results are equivalent if the Golub result is
extended to accomodate complex coefficients and we choose \(c=e_j\). In
this case \(P=I-cc^*\) is a projection matrix which sets the \(j^{th}\)
row and column to zero in \(PAP\), so we can see that the eigenvalues
\(x_k\) are equivalent to the \(x_{jk}\) in the Tao result, and
therefore the \(d_j\) are equivalent to the \(q_{ij}\). Thus the
eigenvector element magnitudes can be constructed by iterating through
all values of \(c=e_j\).

    \hypertarget{arbitrary-orthonormal-basis}{%
\section{Arbitrary Orthonormal
Basis}\label{arbitrary-orthonormal-basis}}

The procedure above is equivalent to choosing \(I\) as the orthonormal
basis for a set of constraints, or equivalently, to form the projection
matrix. We may choose any orthonormal basis \(C^*C=I\) with columns
\(c_j\) forming projection matrix \(P_j=I-c_jc_j^*\). If \(A=QWQ^*\),
\(S=CQ\) has elements \(s_{ij}\), and \(y_{jk}\) denote the nonzero
eigenvalues of \(P_jAP_j\), then
\[ s_{ij}^*s_{ij} = \frac{\prod_{k=1;k\neq i}^{n}{w_i-w_k}}{\prod_{k=1}^{n-1}{w_i-y_{jk}}} \]

    \hypertarget{numerical-experiments}{%
\section{Numerical Experiments}\label{numerical-experiments}}

    \begin{tcolorbox}[breakable, size=fbox, boxrule=1pt, pad at break*=1mm,colback=cellbackground, colframe=cellborder]
\prompt{In}{incolor}{1}{\boxspacing}
\begin{Verbatim}[commandchars=\\\{\}]
\PY{k+kn}{import} \PY{n+nn}{numpy} \PY{k}{as} \PY{n+nn}{np}
\end{Verbatim}
\end{tcolorbox}

    \begin{tcolorbox}[breakable, size=fbox, boxrule=1pt, pad at break*=1mm,colback=cellbackground, colframe=cellborder]
\prompt{In}{incolor}{2}{\boxspacing}
\begin{Verbatim}[commandchars=\\\{\}]
\PY{n}{n} \PY{o}{=} \PY{l+m+mi}{100}
\PY{n}{eps} \PY{o}{=} \PY{l+m+mf}{1e\PYZhy{}10}
\end{Verbatim}
\end{tcolorbox}

    Define a random hermitian matrix and compute the eigenvectors/values

    \begin{tcolorbox}[breakable, size=fbox, boxrule=1pt, pad at break*=1mm,colback=cellbackground, colframe=cellborder]
\prompt{In}{incolor}{3}{\boxspacing}
\begin{Verbatim}[commandchars=\\\{\}]
\PY{n}{A} \PY{o}{=} \PY{n}{np}\PY{o}{.}\PY{n}{random}\PY{o}{.}\PY{n}{random}\PY{p}{(}\PY{p}{(}\PY{n}{n}\PY{p}{,}\PY{n}{n}\PY{p}{)}\PY{p}{)}\PY{o}{+}\PY{n}{np}\PY{o}{.}\PY{n}{random}\PY{o}{.}\PY{n}{random}\PY{p}{(}\PY{p}{(}\PY{n}{n}\PY{p}{,}\PY{n}{n}\PY{p}{)}\PY{p}{)}\PY{o}{*}\PY{l+m+mi}{1}\PY{n}{j}
\PY{n}{A} \PY{o}{=} \PY{n}{A}\PY{o}{+}\PY{n}{A}\PY{o}{.}\PY{n}{conj}\PY{p}{(}\PY{p}{)}\PY{o}{.}\PY{n}{T}
\PY{n}{w}\PY{p}{,}\PY{n}{Q} \PY{o}{=} \PY{n}{np}\PY{o}{.}\PY{n}{linalg}\PY{o}{.}\PY{n}{eigh}\PY{p}{(}\PY{n}{A}\PY{p}{)}
\end{Verbatim}
\end{tcolorbox}

    Define the eigenvector computation function.

    \begin{tcolorbox}[breakable, size=fbox, boxrule=1pt, pad at break*=1mm,colback=cellbackground, colframe=cellborder]
\prompt{In}{incolor}{4}{\boxspacing}
\begin{Verbatim}[commandchars=\\\{\}]
\PY{k}{def} \PY{n+nf}{fR}\PY{p}{(}\PY{n}{W}\PY{p}{,}\PY{n}{w}\PY{p}{)}\PY{p}{:}
    \PY{l+s+sd}{\PYZsq{}\PYZsq{}\PYZsq{}}
\PY{l+s+sd}{    W: eigenvalues of submatrices (n,n\PYZhy{}1)}
\PY{l+s+sd}{    w: eigenvalues of matrix (n,)}
\PY{l+s+sd}{    output: matrix of eigenbasis squared magnitudes (n,n)}
\PY{l+s+sd}{    \PYZsq{}\PYZsq{}\PYZsq{}}
    \PY{k}{return} \PY{n}{np}\PY{o}{.}\PY{n}{prod}\PY{p}{(}
        \PY{n}{w}\PY{p}{[}\PY{p}{:}\PY{p}{,}\PY{n}{np}\PY{o}{.}\PY{n}{newaxis}\PY{p}{,}\PY{n}{np}\PY{o}{.}\PY{n}{newaxis}\PY{p}{]}\PY{o}{\PYZhy{}}\PY{n}{W}\PY{p}{[}\PY{n}{np}\PY{o}{.}\PY{n}{newaxis}\PY{p}{,}\PY{p}{:}\PY{p}{,}\PY{p}{:}\PY{p}{]}\PY{p}{,}
        \PY{n}{axis}\PY{o}{=}\PY{l+m+mi}{2}
    \PY{p}{)}\PY{o}{/}\PY{n}{np}\PY{o}{.}\PY{n}{prod}\PY{p}{(}
        \PY{n}{w}\PY{p}{[}\PY{p}{:}\PY{p}{,}\PY{n}{np}\PY{o}{.}\PY{n}{newaxis}\PY{p}{]}\PY{o}{\PYZhy{}}\PY{n}{w}\PY{p}{[}\PY{n}{np}\PY{o}{.}\PY{n}{newaxis}\PY{p}{,}\PY{p}{:}\PY{p}{]}\PY{o}{+}\PY{n}{np}\PY{o}{.}\PY{n}{eye}\PY{p}{(}\PY{n}{n}\PY{p}{)}\PY{p}{,}
        \PY{n}{axis}\PY{o}{=}\PY{l+m+mi}{1}
    \PY{p}{)}\PY{p}{[}\PY{p}{:}\PY{p}{,}\PY{n}{np}\PY{o}{.}\PY{n}{newaxis}\PY{p}{]}
\end{Verbatim}
\end{tcolorbox}

    \hypertarget{test-of-taos-method}{%
\subsection{Test of Tao's Method}\label{test-of-taos-method}}

    \begin{tcolorbox}[breakable, size=fbox, boxrule=1pt, pad at break*=1mm,colback=cellbackground, colframe=cellborder]
\prompt{In}{incolor}{5}{\boxspacing}
\begin{Verbatim}[commandchars=\\\{\}]
\PY{n}{W} \PY{o}{=} \PY{n}{np}\PY{o}{.}\PY{n}{asarray}\PY{p}{(}\PY{p}{[}
    \PY{n}{np}\PY{o}{.}\PY{n}{linalg}\PY{o}{.}\PY{n}{eigvalsh}\PY{p}{(}
        \PY{n}{A}\PY{p}{[}\PY{p}{[}\PY{n}{j} \PY{k}{for} \PY{n}{j} \PY{o+ow}{in} \PY{n+nb}{range}\PY{p}{(}\PY{n}{n}\PY{p}{)} \PY{k}{if} \PY{n}{j}\PY{o}{!=}\PY{n}{i}\PY{p}{]}\PY{p}{,}\PY{p}{:}\PY{p}{]}\PY{p}{[}\PY{p}{:}\PY{p}{,}\PY{p}{[}\PY{n}{j} \PY{k}{for} \PY{n}{j} \PY{o+ow}{in} \PY{n+nb}{range}\PY{p}{(}\PY{n}{n}\PY{p}{)} \PY{k}{if} \PY{n}{j}\PY{o}{!=}\PY{n}{i}\PY{p}{]}\PY{p}{]}
    \PY{p}{)} \PY{k}{for} \PY{n}{i} \PY{o+ow}{in} \PY{n+nb}{range}\PY{p}{(}\PY{n}{n}\PY{p}{)}
\PY{p}{]}\PY{p}{)}
\PY{n}{R} \PY{o}{=} \PY{n}{fR}\PY{p}{(}\PY{n}{W}\PY{p}{,}\PY{n}{w}\PY{p}{)}
\PY{p}{(}\PY{p}{(}\PY{n}{R}\PY{o}{\PYZhy{}}\PY{p}{(}\PY{n}{Q}\PY{o}{*}\PY{n}{Q}\PY{o}{.}\PY{n}{conj}\PY{p}{(}\PY{p}{)}\PY{p}{)}\PY{o}{.}\PY{n}{T}\PY{p}{)}\PY{o}{\PYZlt{}}\PY{n}{eps}\PY{p}{)}\PY{o}{.}\PY{n}{all}\PY{p}{(}\PY{p}{)}
\end{Verbatim}
\end{tcolorbox}

            \begin{tcolorbox}[breakable, size=fbox, boxrule=.5pt, pad at break*=1mm, opacityfill=0]
\prompt{Out}{outcolor}{5}{\boxspacing}
\begin{Verbatim}[commandchars=\\\{\}]
True
\end{Verbatim}
\end{tcolorbox}
        
    \hypertarget{test-of-equivalence}{%
\subsection{Test of Equivalence}\label{test-of-equivalence}}

    \begin{tcolorbox}[breakable, size=fbox, boxrule=1pt, pad at break*=1mm,colback=cellbackground, colframe=cellborder]
\prompt{In}{incolor}{6}{\boxspacing}
\begin{Verbatim}[commandchars=\\\{\}]
\PY{n}{C} \PY{o}{=} \PY{n}{np}\PY{o}{.}\PY{n}{eye}\PY{p}{(}\PY{n}{n}\PY{p}{)}
\PY{n}{P} \PY{o}{=} \PY{p}{[}\PY{n}{np}\PY{o}{.}\PY{n}{eye}\PY{p}{(}\PY{n}{n}\PY{p}{)}\PY{o}{\PYZhy{}}\PY{n}{np}\PY{o}{.}\PY{n}{outer}\PY{p}{(}\PY{n}{c}\PY{p}{,}\PY{n}{c}\PY{o}{.}\PY{n}{conj}\PY{p}{(}\PY{p}{)}\PY{p}{)} \PY{k}{for} \PY{n}{c} \PY{o+ow}{in} \PY{n}{C}\PY{p}{]}
\PY{n}{W} \PY{o}{=} \PY{n}{np}\PY{o}{.}\PY{n}{asarray}\PY{p}{(}\PY{p}{[}
    \PY{n}{v}\PY{p}{[}\PY{n}{np}\PY{o}{.}\PY{n}{argsort}\PY{p}{(}\PY{n}{np}\PY{o}{.}\PY{n}{abs}\PY{p}{(}\PY{n}{v}\PY{p}{)}\PY{p}{)}\PY{p}{[}\PY{l+m+mi}{1}\PY{p}{:}\PY{p}{]}\PY{p}{]} 
    \PY{k}{for} \PY{n}{v} \PY{o+ow}{in} \PY{p}{[}
        \PY{n}{np}\PY{o}{.}\PY{n}{linalg}\PY{o}{.}\PY{n}{eigvalsh}\PY{p}{(}\PY{n}{np}\PY{o}{.}\PY{n}{linalg}\PY{o}{.}\PY{n}{multi\PYZus{}dot}\PY{p}{(}\PY{p}{[}\PY{n}{p}\PY{p}{,}\PY{n}{A}\PY{p}{,}\PY{n}{p}\PY{p}{]}\PY{p}{)}\PY{p}{)}
        \PY{k}{for} \PY{n}{p} \PY{o+ow}{in} \PY{n}{P}
    \PY{p}{]}
\PY{p}{]}\PY{p}{)}
\PY{n}{R} \PY{o}{=} \PY{n}{fR}\PY{p}{(}\PY{n}{W}\PY{p}{,}\PY{n}{w}\PY{p}{)}
\PY{p}{(}\PY{p}{(}\PY{n}{R}\PY{o}{\PYZhy{}}\PY{p}{(}\PY{n}{Q}\PY{o}{*}\PY{n}{Q}\PY{o}{.}\PY{n}{conj}\PY{p}{(}\PY{p}{)}\PY{p}{)}\PY{o}{.}\PY{n}{T}\PY{p}{)}\PY{o}{\PYZlt{}}\PY{n}{eps}\PY{p}{)}\PY{o}{.}\PY{n}{all}\PY{p}{(}\PY{p}{)}
\end{Verbatim}
\end{tcolorbox}

            \begin{tcolorbox}[breakable, size=fbox, boxrule=.5pt, pad at break*=1mm, opacityfill=0]
\prompt{Out}{outcolor}{6}{\boxspacing}
\begin{Verbatim}[commandchars=\\\{\}]
True
\end{Verbatim}
\end{tcolorbox}
        
    \hypertarget{test-of-arbitrary-orthonormal-basis}{%
\subsection{Test of Arbitrary Orthonormal
Basis}\label{test-of-arbitrary-orthonormal-basis}}

    \begin{tcolorbox}[breakable, size=fbox, boxrule=1pt, pad at break*=1mm,colback=cellbackground, colframe=cellborder]
\prompt{In}{incolor}{7}{\boxspacing}
\begin{Verbatim}[commandchars=\\\{\}]
\PY{k+kn}{from} \PY{n+nn}{scipy}\PY{n+nn}{.}\PY{n+nn}{stats} \PY{k}{import} \PY{n}{ortho\PYZus{}group}
\PY{n}{C} \PY{o}{=} \PY{n}{ortho\PYZus{}group}\PY{o}{.}\PY{n}{rvs}\PY{p}{(}\PY{n}{n}\PY{p}{)}
\PY{n}{P} \PY{o}{=} \PY{p}{[}\PY{n}{np}\PY{o}{.}\PY{n}{eye}\PY{p}{(}\PY{n}{n}\PY{p}{)}\PY{o}{\PYZhy{}}\PY{n}{np}\PY{o}{.}\PY{n}{outer}\PY{p}{(}\PY{n}{c}\PY{p}{,}\PY{n}{c}\PY{o}{.}\PY{n}{conj}\PY{p}{(}\PY{p}{)}\PY{p}{)} \PY{k}{for} \PY{n}{c} \PY{o+ow}{in} \PY{n}{C}\PY{p}{]}
\PY{n}{W} \PY{o}{=} \PY{n}{np}\PY{o}{.}\PY{n}{asarray}\PY{p}{(}\PY{p}{[}
    \PY{n}{v}\PY{p}{[}\PY{n}{np}\PY{o}{.}\PY{n}{argsort}\PY{p}{(}\PY{n}{np}\PY{o}{.}\PY{n}{abs}\PY{p}{(}\PY{n}{v}\PY{p}{)}\PY{p}{)}\PY{p}{[}\PY{l+m+mi}{1}\PY{p}{:}\PY{p}{]}\PY{p}{]} 
    \PY{k}{for} \PY{n}{v} \PY{o+ow}{in} \PY{p}{[}
        \PY{n}{np}\PY{o}{.}\PY{n}{linalg}\PY{o}{.}\PY{n}{eigvalsh}\PY{p}{(}\PY{n}{np}\PY{o}{.}\PY{n}{linalg}\PY{o}{.}\PY{n}{multi\PYZus{}dot}\PY{p}{(}\PY{p}{[}\PY{n}{p}\PY{p}{,}\PY{n}{A}\PY{p}{,}\PY{n}{p}\PY{p}{]}\PY{p}{)}\PY{p}{)}
        \PY{k}{for} \PY{n}{p} \PY{o+ow}{in} \PY{n}{P}
    \PY{p}{]}
\PY{p}{]}\PY{p}{)}
\PY{n}{R} \PY{o}{=} \PY{n}{fR}\PY{p}{(}\PY{n}{W}\PY{p}{,}\PY{n}{w}\PY{p}{)}
\PY{n}{S} \PY{o}{=} \PY{n}{np}\PY{o}{.}\PY{n}{dot}\PY{p}{(}\PY{n}{C}\PY{p}{,}\PY{n}{Q}\PY{p}{)}
\PY{p}{(}\PY{p}{(}\PY{n}{R}\PY{o}{\PYZhy{}}\PY{p}{(}\PY{n}{S}\PY{o}{*}\PY{n}{S}\PY{o}{.}\PY{n}{conj}\PY{p}{(}\PY{p}{)}\PY{p}{)}\PY{o}{.}\PY{n}{T}\PY{p}{)}\PY{o}{\PYZlt{}}\PY{n}{eps}\PY{p}{)}\PY{o}{.}\PY{n}{all}\PY{p}{(}\PY{p}{)}
\end{Verbatim}
\end{tcolorbox}

            \begin{tcolorbox}[breakable, size=fbox, boxrule=.5pt, pad at break*=1mm, opacityfill=0]
\prompt{Out}{outcolor}{7}{\boxspacing}
\begin{Verbatim}[commandchars=\\\{\}]
True
\end{Verbatim}
\end{tcolorbox}
        
    \hypertarget{discussion}{%
\section{Discussion}\label{discussion}}

We have shown equivalence of the methods of Denton-Tao and Golub for
computing eigenvector element magnitudes from eigenvalues resulting from
row-column eliminations, and shown a generalization of these
relationships to projections composed from an arbitrary orthonormal
basis. We have produced program code to test these methods, and provided
evidence by computing results of randomized inputs.


\bibliographystyle{unsrt}
\bibliography{GolubTao}

\begin{thebibliography}{1}

\bibitem{denton2019eigenvectors}
Peter~B. Denton, Stephen~J. Parke, Terence Tao, and Xining Zhang.
\newblock Eigenvectors from eigenvalues, 2019.

\bibitem{10.2307/2028604}
Gene~H. Golub.
\newblock Some modified matrix eigenvalue problems.
\newblock {\em SIAM Review}, 15(2):318--334, 1973.

\end{thebibliography}

\end{document}